\documentclass[11pt]{article}
\def\be{\begin{equation}}
\def\ee{\end{equation}}
\newcommand{\ff}[1]{{\mbox{\boldmath $#1$}}}

\def\x{\ff{x}}
\def\pk{\ff{p}}
\def\vem{\ff{\varepsilon}}
\newcommand{\vcode}[2]{
\begin{figure} \begin{center}
\begin{minipage}{0.8\textwidth}
\rule{\linewidth}{1pt}
#2
\rule{\linewidth}{1pt}
\caption{#1}
\end{minipage}
\end{center} \end{figure}
}    

\begin{document}
\title{A Framework for Self-Tuning Optimization Algorithm }
\author{Xin-She Yang$^1$, Suash Deb,  Martin Loomes$^1$, Mehmet Karamanoglu$^1$ \\ \\
{\small 1) School of Science and Technology, Middlesex University,  London NW4 4BT, UK. }\\
{\small 2) Cambridge Institute of Technology,  Cambridge Village, Tatisilwai,  India.
}}

\date{}

\maketitle

\begin{abstract}
The performance of any algorithm will largely depend on the setting of its
algorithm-dependent parameters. The optimal setting should allow the algorithm
to achieve the best performance for solving a range of optimization problems.
However, such parameter-tuning itself is a tough optimization problem. In this paper,
we present a framework for self-tuning algorithms so that an algorithm to be tuned
can be used to tune the algorithm itself. Using the firefly algorithm as an example,
we show that this framework works well. It is also found that different parameters
may have different sensitivities, and thus require different degrees of tuning.
Parameters with high sensitivities require fine-tuning to achieve optimality. \\


\end{abstract}

{\bf Citation Details:} X. S. Yang, S. Deb, M. Loomes, M. Karamanoglu, A framework for self-tuning 
optimization algorithms, {\it Neural Computing and Applications}, Vol. 23, No. 7-8, pp. 2051-2057 (2013).

\section{Introduction}

Optimization is paramount in many applications such as engineering and industrial designs.
Obviously, the aims of optimization can be anything -- to minimize the energy consumption,
to maximize the profit, output, performance and efficiency \cite{Yang2008book,KozielYang,Yang2010book,YangAPSO}.
As most real-world applications are often highly nonlinear, it requires sophisticated
optimization tools to tackle. There are many algorithms that use swarm intelligence to solve
optimization problems, and algorithms such as particle swarm optimization,
cuckoo search and firefly algorithm have received a lot of interests.
These nature-inspired algorithms have been proved very efficient.

Metaheuristic algorithms are often nature-inspired, and they are now among the most widely used algorithms for optimization.
They have many advantages over conventional algorithms \cite{KozielYang,Yang2010book,Gandomi}.
Metaheuristic algorithms are very diverse, including genetic algorithms,
simulated annealing, differential evolution, ant and bee algorithms,
bat algorithm, particle swarm optimization, harmony search, firefly algorithm, cuckoo search and others
\cite{Kennedy,YangFA,YangBA2012,Gandomi2}.

Since all algorithms have algorithm-dependent parameters, the performance of an algorithm
largely depends on the values or setting of these parameters. Ideally, there
should be a good way to tune these parameters so that the performance of the algorithm
can be optimal in the sense that the algorithm can find the optimal solution of
a problem using the minimal number of iterations and with the highest accuracy.
However, such tuning of algorithm-dependent parameters is itself a very tough
optimization problem. In essence, it is a hyper-optimization problem, that is
the optimization of optimization. In fact, how to find the best parameter setting
of an algorithm is still an open problem.

There are studies on parameter tuning. For example, Eiben provided a comprehensive
summary of existing studies \cite{Eiben}. However, these studies are still very preliminary.
There is no method of self-tuning in algorithms. Therefore, the main objective
of this paper is to provide a framework for self-tuning algorithms so that
an algorithm can be used to tune its own parameters automatically. As far as we are concerned,
this is the first of its kind in parameter tuning. The paper is thus organized as follows:
Section 2 first analyzes the essence of parameter tuning and Section 3 provides
a framework for automatic parameter tuning.
Section 4 uses the firefly algorithm to show how the self-tuning framework works.
Then, Section 5 presents a case study of a gearbox design problem to further
test the tuning procedure. Finally, we draw conclusions briefly in Section 6.

\section{Algorithm Analysis and Parameter Tuning}

An optimization algorithm is essentially an iterative procedure, starting
with some initial guess point/solution with an aim to reach a better
solution or ideally the optimal solution to a problem of interest.
This process of search for optimality is generic, though
the details of the process can vary from algorithm to algorithm.
Traditional algorithms such as Newton-Raphson methods use a deterministic
trajectory-based method, while modern nature-inspired algorithms often
are population-based algorithms that use multiple agents. In essence,
these multiple agents form an iterative, dynamic system which should
have some attractors or stable states. On the other hand, the same
system can be considered as a set of Markov chains so that they will
converge towards some stable probability distribution.

\subsection{An Optimization Algorithm}

Whatever the perspective may be,  the aim of such an
iterative process is to let the evolve system and converge into
some stable optimality. In this case, it has strong similarity to
a self-organizing system. Such an iterative, self-organizing system
can evolve, according to a set of rules or mathematical equations.
As a result,  such a complex system can interact and self-organize into certain
converged states, showing some emergent characteristics of self-organization.
In this sense, the proper design of an efficient optimization algorithm
is equivalent to finding efficient ways to mimic
the evolution of a self-organizing system \cite{Ashby,Keller}.

From a mathematical point of view, an algorithm $A$ tends to generate
a new and better solution $\x^{t+1}$ to a given problem from the current
solution $\x^t$ at iteration or time $t$. For example, the Newton-Raphson method
to find the optimal solution of $f(\x)$ is equivalent to finding the critical points
or roots of $f'(\x)=0$ in a $d$-dimensional space. That is,
\be \x^{t+1} = \x^t -\frac{f'(\x^t)}{f''(\x^t)}=A(\x^t). \ee
Obviously, the convergence rate may become very slow near the optimal point
where $f'(x) \rightarrow 0$. In general, this Newton-Raphson method has
a quadratic convergence rate \cite{Yang2008}. Sometimes, the true convergence rate may not
be as quick as it should be, it may have non-quadratic convergence property.
A way to improve the convergence in this case is to modify the above formula
slightly by introducing a parameter $p$ so that
\be \x^{t+1} = \x^t -p \frac{f'(\x^t)}{f''(\x^t)}. \ee
If the optimal solution, i.e., the fixed point of the iterations \cite{Suli}, is
$\x_*$, then we can take $p$ as
\be p=\frac{1}{1-A'(\x_*)}. \ee
The above iterative equation can be written as
\be \x^{t+1}=A(\x^t,p). \ee
It is worth pointing out that the optimal convergence of Newton-Raphson's method
leads to an optimal parameter setting $p$ which depends on the iterative
formula and the optimality $\x_*$ of the objective $f(\x)$ to be optimized.

This above formula is valid for a deterministic method; however, in modern
metaheuristic algorithms, randomization is often used in an algorithm, and
in many cases, randomization appears in the form of a set of $m$ random variables $\vem=(\varepsilon_1,...,\varepsilon_m)$
in an algorithm. For example, in simulated annealing, there is one random variable,
while in particle swarm optimization \cite{Kennedy}, there are two random variables.
In addition, there are often a set of $k$ parameters in an algorithm. For example,
in particle swarm optimization, there are 4 parameters (two learning parameters,
one inertia weight, and the population size). In general, we can have a vector
of parameters $\pk=(p_1, ..., p_k)$. Mathematically speaking, we can write an
algorithm with $k$ parameters and $m$ random variables as
\be \x^{t+1}=\ff{A}\Big(\x^t, \ff{p}(t), \ff{\varepsilon}(t)\Big), \label{Equ-Alg} \ee
where $\ff{A}$ is a nonlinear mapping from a given solution (a $d$-dimensional
vector $\x^t$) to a new solution vector $\x^{t+1}$.

\subsection{Type of Optimality}

Representation (\ref{Equ-Alg}) gives rise to two types of optimality:
optimality of a problem and optimality of an algorithm. For an optimization problem such as
$\min f(\x)$, there is a global optimal solution whatever the algorithmic tool we may use
to find this optimality. This is the optimality for the optimization problem.  On the other hand,
for a given problem $\Phi$ with an objective function $f(\x)$,
there are many algorithms that can solve it.
Some algorithms may require less computational effort than others. There may be
the best algorithm with the least computing cost, though this may not be unique.
However, this is not our concern here. Once we have chosen an algorithm $A$
to solve a problem $\Phi$, there is an optimal parameter setting for this algorithm
so that it can achieve the best performance. This optimality depends on both the algorithm
itself and the problem it solves. In the rest of this paper, we will focus on this type
of optimality.

That is,  the optimality to be achieved is
\be \textrm{Maximize the performance of } \xi=\; A(\Phi, \pk, \vem), \ee
for a given problem $\Phi$ and a chosen algorithm $A(., \pk, \vem)$.
We will denote this optimality as $\xi_*=A_*(\Phi,\pk_*)=\xi(\Phi, \pk_*)$
where $\pk_*$ is the optimal parameter setting for this algorithm so that
its performance is the best. Here, we have used a fact that $\vem$ is a random vector can be drawn from some
known probability distributions, thus the randomness vector should not be related to
the algorithm optimality.

It is worth pointing out that there is another potential optimality. That is,
for a given problem, a chosen algorithm with the best parameter setting $\pk_*$,
we can still use different random numbers drawn from various probability distributions
and even chaotic maps, so that the performance can achieve even better performance.
Strictly speaking, if an algorithm $A(.,.,\vem)$  has a random vector $\vem$ that is
drawn from a uniform distribution $\vem_1 \sim U(0,1)$ or from a Gaussian $\vem_2 \sim N(0,1)$,
it becomes two algorithms $A_1=A(.,.,\vem_1)$ and $A_2=A(.,.,\vem_2)$. Technically speaking,
we should treat them as different algorithms.
Since our emphasis here is about parameter tuning so as to find the optimal setting of parameters,
we will omit effect of the randomness vector, and thus focus on
\be \textrm{ Maximize } \xi=A(\Phi, \pk). \ee
In essence, tuning algorithm involves in tuning its algorithm-dependent parameters. Therefore,
parameter tuning is equivalent to algorithm tuning in the present context.

\subsection{Parameter Tuning }

In order to tune $A(\Phi, \pk)$ so as to achieve its best performance, a parameter-tuning tool,
i.e., a tuner, is needed. Like tuning a high-precision machinery, sophisticated tools are required.
For tuning parameters in an algorithm, what tool can we use? One way is to use
a better, existing tool (say, algorithm $B$) to tune an algorithm $A$. Now the question
may become: how do you know $B$ is better? Is $B$ well-tuned? If yes, how do you tune
$B$ in the first place? Naively, if we say, we use another tool (say, algorithm $C$)
to tune $B$. Now again the question becomes how algorithm $C$ has been tuned?
This can go on and on, until the end of a long chain, say, algorithm $Q$.
In the end, we need some tool/algorithm to tune this $Q$, which again come back to
the original question: how to tune an algorithm $A$ so that it can perform best.

It is worth pointing out that even if we have good tools to tune an algorithm, the best
parameter setting and thus performance all depend on the performance measures used
in the tuning. Ideally, the parameters should be robust enough to minor parameter
changes, random seeds, and even problem instance \cite{Eiben}. However, in practice,
they may not be achievable. According to Eiben \cite{Eiben}, parameter tuning can be
divided into iterative and non-iterative tuners, single-stage and multi-stage tuners.
The meaning of these terminologies is self-explanatory. In terms of the actual tuning
methods, existing methods include sampling methods, screening methods,
model-based methods, and metaheuristic methods. Their success and effectiveness can vary, and thus there are no
well-established methods for universal parameter tuning.

\section{Framework for Self-Tuning Algorithms}

\subsection{Hyper-optimization}

From our earlier observations and discussions, it is clear that parameter tuning is the process
of optimizing the optimization algorithm, therefore, it is a hyper-optimization problem. In essence,
a tuner is a meta-optimization tool for tuning algorithms.

For a standard unconstrained optimization problem, the aim is to find the global minimum $f_*$
of a function $f(\x)$ in a $d$-dimensional space. That is,
\be \textrm{Minimize } \; f(\x), \quad \x=(x_1, x_2, ..., x_d). \ee
Once we choose an algorithm $A$ to solve this optimization problem, the algorithm will find
a minimum solution $f_{\min}$ which may be close to the true global minimum $f_*$.
For a given tolerance $\delta$, this may requires $t_{\delta}$ iterations to achieve $|f_{\min}-f_*| \le \delta$.
Obviously, the actual $t_{\delta}$ will largely depend on both the problem objective
$f(\x)$ and the parameters $\pk$ of the algorithm used.

The main aim of algorithm-tuning is to find the best parameter setting $\pk_*$ so that
the computational cost or the number of iterations $t_{\delta}$ is the minimum.  Thus, parameter tuning
as a hyper-optimization problem can be written as
\be \textrm{Minimize } \; t_{\delta}=A(f(\x), \pk), \ee
whose optimality is $\pk_*$.

Ideally, the parameter vector $\pk_*$ should be sufficiently robust. For different types of problems,
any slight variation in $\pk_*$ should not affect the performance of $A$ much, which means
that $\pk_*$ should lie in a flat range, rather than at a sharp peak in the
parameter landscape.

\subsection{Multi-Objective View}

If we look the algorithm tuning process from a different perspective, it is possible to
construct it as a multi-objective optimization problem with two objectives: one
objective $f(\x)$ for the problem $\Phi$ and one objective $t_{\delta}$ for
the algorithm. That is
\be \textrm{Minimize } f(\x) \textrm{ and Minimize }  t_{\delta}=A(f(\x),\pk), \label{two-min-equ} \ee
where $t_{\delta}$ is the (average) number of iterations needed to achieve
 a given tolerance $\delta$ so that the found minimum $f_{\min}$
 is close enough to the true global minimum $f_*$, satisfying $|f_{\min}-f_*| \le \delta$.

This means that for a given tolerance $\delta$, there will be a set of best parameter settings
with a minimum $t_{\delta}$. As a result, the bi-objectives will form a Pareto front.
In principle, this  bi-objective optimization problem (\ref{two-min-equ}) can be solved by
any methods that are suitable for multiobjective optimization. But as $\delta$ is usually given,
a natural way to solve this problem is to use the so-called $\epsilon$-constraint or $\delta$-constraint
methods. The naming may be dependent on the notations; however, we will use $\delta$-constraints.

For a given $\delta \ge 0$, we change one of the objectives (i.e., $f(\x)$) into a constraint,
and thus the above problem (\ref{two-min-equ}) becomes a single-objective optimization problem
with a constraint. That is
\be \textrm{Minimize }\; t_{\delta}=A(f(\x), \pk),  \label{para-equ-single} \ee
subject to
\be f(\x) \le \delta.  \ee
In the rest of this paper, we will set $\delta=10^{-5}$.

The important thing is that we still need an algorithm to solve this optimization problem.
However, the main difference from a common single objective problem is that
the present problem contains an algorithm $A$. Ideally, an algorithm should be
independent of the problem, which treats the objective to be solved
as a black box. Thus we have $A(.,\pk, \vem)$, however, in reality, an algorithm
will be used to solve a particular problem $\Phi$ with an objective $f(\x)$.
Therefore, both notations $A(., \pk)$ and $A(f(\x),\pk)$ will be used in this paper.

\subsection{Self-Tuning Framework }
In principle, we can solve (\ref{para-equ-single}) by any efficient
or well-tuned algorithm. Now a natural question is: Can we solve this algorithm-tuning problem
by the algorithm $A$ itself? There is no reason we cannot. In fact, if we solve (\ref{para-equ-single}) by using
$A$, we have a self-tuning algorithm. That is, the algorithm automatically tunes itself for a given problem objective
to be optimized. This essentially provides a framework for a self-tuning algorithm as shown in Fig.~\ref{Fig-framework}.

\vcode{A Framework for a Self-Tuning Algorithm. \label{Fig-framework}}
{Implement an algorithm $A(.,\ff{p}, \ff{\varepsilon})$  with $\pk=[p_1, ..., p_K], \vem=[\varepsilon_1, ..., \varepsilon_m]$; \\
Define a tolerance (e.g., $\delta=10^{-5}$); \\
\indent $\quad$ Algorithm objective $t_{\delta}(f(\x), \pk, \vem)$; \\
\indent $\qquad$ Problem objective function $f(\x)$; \\
\indent $\qquad$ Find the optimality solution $f_{\min}$ within $\delta$; \\
\indent $\qquad $ Output the number of iterations $t_{\delta}$ needed to find $f_{\min}$; \\
\indent $\quad$ Solve $\min t_{\delta}(f(\x), \pk)$ using $A(., \pk, \vem)$ to get the best parameters;  \\
Output the tuned algorithm with the best parameter setting $\pk_*$. \\
}

This framework is generic in the sense that any algorithm can be tuned this way, and
any problem can be solved within this framework. This essentially achieves two goals simultaneously:
parameter tuning and optimality finding.

In the rest of this paper, we will use firefly algorithm (FA) as a case study to self-tune FA
for a set of function optimization problems.

\section{Self-Tuning Firefly Algorithm}

\subsection{Firefly Algorithm}

Firefly Algorithm (FA) was developed by Xin-She Yang in 2008 \cite{Yang2008book,YangFA,YangFA2},
which was based on the flashing patterns and behaviour of tropical fireflies.
In essence, FA uses the following three idealized rules:

\begin{itemize}
\item Fireflies are unisex so that one firefly will be attracted to other fireflies regardless of their sex.

\item The attractiveness is proportional to the brightness and they both decrease as their distance increases.
Thus for any two flashing fireflies, the less brighter one will move towards the brighter one.
If there is no brighter one than a particular firefly, it will move randomly.

\item The brightness of a firefly is determined by the landscape of the objective function.
\end{itemize}

As a firefly's attractiveness is proportional to the light intensity seen by adjacent fireflies, we can now define the variation of attractiveness $\beta$ with the distance $r$ by
    \be \beta = \beta_0 e^{-\gamma r^2}, \ee
where $\beta_0$ is the attractiveness at $r=0$.

The movement of a firefly $i$ is attracted to another more attractive (brighter) firefly $j$ is determined by
\be    \x_i^{t+1} =\x_i^t + \beta_0 e^{-\gamma r^2_{ij} } (\x_j^t-\x_i^t) + \alpha \; \ff{\epsilon}_i^t, \ee
where the second term is due to the attraction. The third term is randomization with $\alpha$ being the randomization parameter, and $\ff{\epsilon}_i^t$ is a vector of random numbers drawn from a Gaussian distribution at time $t$. Other studies also use
the randomization in terms of $\ff{\epsilon}_i^t$ that can easily be extended to other distributions such as L\'evy flights \cite{YangFA,YangFA2}.

For simplicity for parameter tuning, we assume that $\beta_0=1$, and therefore the two parameters to be tuned are:
$\gamma>0$ and $\alpha>0$. It is worth pointing out that $\gamma$ controls the scaling, while $\alpha$ controls
the randomness. For this algorithm to convergence properly, randomness should be gradually reduced, and one way to
achieve such randomness reduction is to use
\be \alpha=\alpha_0 \theta^t, \quad \theta \in (0, 1), \ee
where $t$ is the index of iterations/generations.
Here $\alpha_0$ is the initial randomness factor, and we can set $\alpha_0=1$ without losing generality.
Therefore, the two parameters to be tuned become $\gamma$ and $\theta$.

\subsection{Tuning the Firefly Algorithm}

Now we will use the framework outlined earlier in this paper to tune FA for a set of five test functions.
The Ackley function can be written as
\be f_1(\x)=-20 \exp\Big[-\frac{1}{5} \sqrt{\frac{1}{d} \sum_{i=1}^d x_i^2}\Big] - \exp\Big[\frac{1}{d} \sum_{i=1}^d \cos (2 \pi x_i)\Big]
+20 +e, \ee
which has a global minimum $f_*=0$ at $(0,0,...,0)$.

The simplest of De Jong's functions is the so-called sphere function
\be f_2(\x) =\sum_{i=1}^d x_i^2, \quad -5.12 \le x_i \le 5.12, \ee
whose global minimum is obviously $f_*=0$ at $(0,0,...,0)$. This function is unimodal and convex.

Yang's forest function \cite{YangFA2}
\be f_3(\x)=\Big( \sum_{i=1}^d |x_i| \Big) \exp\Big[- \sum_{i=1}^d \sin (x_i^2) \Big],
\;\;\; -2 \pi \le x_i \le 2 \pi, \ee
is highly multimodal and has a global minimum $f_*=0$ at $(0,0,...,0)$.

Rastrigin's function \be f_4(\x) = 10 d + \sum_{i=1}^d \Big[ x_i^2 - 10 \cos (2 \pi x_i) \Big],
\quad -5.12 \le x_i \le 5.12, \ee
whose global minimum is $f_*=0$ at $(0,0,...,0)$. This function is highly multimodal. \\

Zakharov's function \cite{Yang2010book}
\be f_5(\x)=\sum_{i=1}^d x_i^2 +\Big(\frac{1}{2} \sum_{i=1}^d i x_i \Big)^2
+\Big(\frac{1}{2} \sum_{i=1}^d i x_i \Big)^4, \ee
has a global minimum $f_*=0$ at $(0,0,...,0)$.

For each objective function, we run the FA to tune itself 50 times so as to
calculated meaningful  statistics. The population size $n=20$ is used for all the
runs.  The means and standard deviations are summarized
in Table \ref{table-mean} where $d=8$ is used for all functions.

\begin{table}
\begin{center}
\caption{Results of parameter tuning for the firefly algorithm. \label{table-mean}}
\begin{tabular}{|l|l|l|l|l|}
\hline
Function & Mean $t_{\delta} \pm \sigma_{t}$ & Mean $\gamma \pm \sigma_{\gamma}$ & Mean $\theta \pm \sigma_{\theta}$ \\ \hline
$f_1$ & $589.7 \pm 182.1$ & $0.5344 \pm 0.2926$ & $0.9561 \pm 0.0076$ \\ \hline
$f_2$ & $514.4 \pm 178.5$ & $0.5985 \pm 0.2554$ & $0.9540 \pm 0.0072$ \\ \hline
$f_3$ & $958.1 \pm 339.0$ & $1.0229 \pm 0.5762$ & $0.9749 \pm 0.0047$ \\ \hline
$f_4$ & $724.1 \pm 217.6$ & $0.4684 \pm 0.3064$ & $0.9652 \pm 0.0065$ \\ \hline
$f_5$ & $957.2 \pm 563.6$ & $0.8933 \pm 0.4251$ & $0.9742 \pm 0.0052$ \\ \hline
\end{tabular}
\end{center}
\end{table}

From this table, we can see that the variations of $\gamma$ is large, while
$\theta$ has a narrow range. The best settings for parameters are problem-dependent.
These results imply the following:
\begin{itemize}
\item[$\bullet$] The optimal setting of parameters in an algorithm largely depends on the
problem, and there is no unique best setting for all problems.

\item[$\bullet$] The relatively large standard deviation of $\gamma$ means that the actual
setting of $\gamma$ is not important to a given problem, and therefore, there is no need
to fine tune $\gamma$. That is to say, a typical value of $\gamma=1$ should work for most problems.

\item[$\bullet$] Some parameters are more sensitive than others. In the present case,
$\theta$ needs more fine-tuning, due to its smaller standard deviations.
\end{itemize}
These findings confirm the earlier observations in the literature that $\gamma=O(1)$ can
be used for most applications \cite{Yang2008book,YangFA}, while $\alpha$ needs to reduce gradually  in terms
of $\theta$. That is probably why other forms of probability distributions such as L\'evy flights
may lead to better performance then the random numbers drawn from the Gaussian normal distribution \cite{YangFA2}.

\section{Applications}

From the results for the test functions, we know that the tuning of $\gamma$ is not important, while
$\theta$ needs more fine-tuning. Let us see if this conclusion is
true for a real-world application. In the rest of the paper, let us focus on
a gearbox design problem.

The optimal design of a speed reducer or a gearbox is a well-known design benchmark with
seven design variables \cite{Cag,Gandomi},
including the face width ($b$), module of the teeth ($h$),
the number of teeth on pinion ($z$),
the length ($L_1$) of the first shaft between bearing, the length ($L_2$)
of the second shaft between between bearings, the diameter ($d_1$) of the first shaft,
and the diameter ($d_2$) of the second shaft. The main objective is to
minimize the total weight of the speed reducer, subject to 11 constraints such as
bending stress, deflection and various limits on stresses in shafts.
This optimization problem can be written as
\[ f(b,h,z,L_1,L_2,d_1,d_2)=0.7854 b h^2 (3.3333 z^2 +14.9334 z -43.0934) \]
\be -1.508 b (d_1^2 +d_2^2) + 7.4777 (d_1^3+d_2^3) + 0.7854 (L_1 d_1^2 + L_2 d_2^2), \ee
subject to
\be \begin{array}{ll} g_1 =\frac{27}{b h^2 z} -1 \le 0,  &
g_2=\frac{397.5}{b h^2 z^2}-1 \le 0, \\ \\
g_3=\frac{1.93L_1^3}{h z d_1^4}- 1 \le 0, &
g_4=\frac{1.93 L_2^3}{h z d_2^4} -1 \le 0, \\ \\
g_5=\frac{1}{110 d_1^3} \sqrt{(\frac{745 L_1}{h z})^2 +16.9 \times 10^6} -1 \le 0, \\ \\
g_6=\frac{1}{85 d_2^3} \sqrt{(\frac{745 L_2}{ h z})^2 +157.5 \times 10^6}-1 \le 0, \\ \\
g_7=\frac{h z}{40}-1 \le 0, &
g_8=\frac{5 h}{b}-1 \le 0, \\ \\
g_9=\frac{b}{12 h} -1 \le 0, &
g_{10}=\frac{1.5 d_1 +1.9}{L_1} -1 \le 0, \\ \\
g_{11}=\frac{1.1 d_2 + 1.9}{L_2} -1 \le 0.
\end{array}  \ee
In addition, the simple bounds are $2.6 \le b \le 3.6, 0.7 \le h \le 0.8$,
$17 \le z \le 28, 7.3 \le L_1 \le 8.3$, $7.8 \le L_2 \le 8.3, 2.9 \le d_1 \le 3.9$,
and $5.0 \le d_2 \le 5.5$. $z$ must be integers.

By using the self-tuning framework via the firefly algorithm with $n=20$, the following best solutions have been obtained:
\[ b=3.5, \; h=0.7, \; z=17, \; L_1=7.3, \;L_2=7.8, \]
\be d_1=3.34336445,\; d_2=5.285350625, \quad f{\min}=2993.7495888, \ee
which are better than $f_*=2996.348165$ obtained by others \cite{Cag,Gandomi}.

The best parameters obtained after tuning are $\gamma=1.0279 \pm 0.4937$ and $\theta=0.9812 \pm 0.0071$,
which are indeed consistent with the results in Table~\ref{table-mean}.

\section{Discussion}

Parameter tuning is the process of tuning an algorithm to find the best parameter settings so that
an algorithm can perform the best for a given set of problems. However, such parameter tuning
is a very tough optimization problem. In fact, such hyper-optimization is
the optimization of an optimization algorithm, which requires special care because the
optimality depends on both the algorithm to be tuned and the problem to be solved.
Though it is possible to view this parameter-tuning process as a bi-objective optimization
problem; however, the objectives involve an algorithm and thus this bi-objective problem is
different from the multiobjective problem in the normal sense.

In this paper, we have successfully developed a framework for self-tuning algorithms in the sense
that the algorithm to be tuned is used to tune itself. We have used the firefly algorithm and
a set of test functions to test the proposed self-tuning algorithm framework. Results have shown that
it can indeed work well. We also found that some parameters require fine-tuning, while others do not
need to be tuned carefully. This is because different parameters may have different sensitivities,
and thus may affect the performance of an algorithm in different ways. Only parameters with high
sensitivities need careful tuning.

Though successful, the present framework requires further extensive testing with
a variety of test functions and many different algorithms. It may also be possible to
see how probability distributions can affect the tuned parameters and even the parameter tuning
process. In addition, it can be expected that this present framework is also useful for parameter control,
so a more generalized framework for both parameter tuning and control can be used for
a wide range of applications. Furthermore, our current framework may be extended to
multiobjective problems so that algorithms for multiobjective optimization can be tuned
in a similar way.

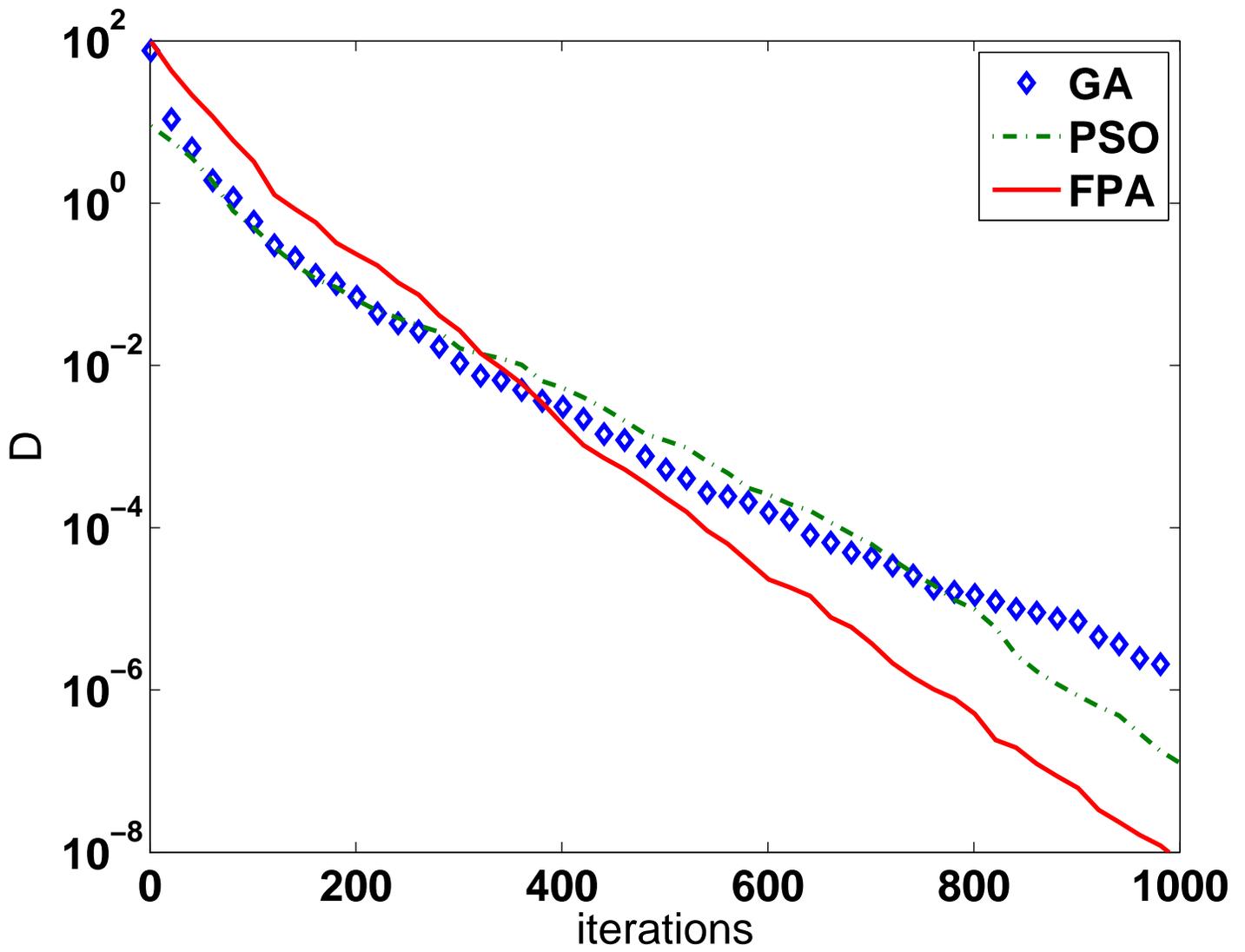
\end{document}